\documentclass{article}
\usepackage{amsmath,amsthm,amsfonts,amssymb,amscd,enumerate}
\theoremstyle{plain}
\newtheorem{theorem}{Theorem}[section]
\newtheorem{proposition}[theorem]{Proposition}

\newtheorem{corollary}[theorem]{Corollary}
\newtheorem{lemma}[theorem]{Lemma}

\theoremstyle{definition}

\newtheorem{example}[theorem]{Example}

\newtheorem*{Remark}{Remark}

        {\end{enumerate}}

\newenvironment{enumeratea'}{\begin{enumerate}[\upshape
(a)$'$]}{\end{enumerate}}

\newcommand{\comment}[1]{}

\numberwithin{equation}{section}
\begin{document}
\title{On Irreducible, Infinite, Non-affine Coxeter Groups}
\author{{Dongwen Qi}\\
{\small The Ohio State University, Department of Mathematics}\\
{\small  231 West 18th Avenue, Columbus, OH 43210, U.S.A.}}
\date{}
\maketitle
\begin{center}
\begin{minipage}{130mm}
\vskip 0.6cm
\begin{center}{\bf Abstract}\end{center}
{The following results are proved:  The center of any finite
index subgroup of an irreducible, infinite, non-affine Coxeter group
is trivial;  Any finite index subgroup of an irreducible,
infinite, non-affine Coxeter group cannot be expressed as a product
of two nontrivial subgroups. These two theorems imply a unique
decomposition theorem for a class of Coxeter groups. We also obtain that
the orbit of each element other than the identity under the conjugation
action in an irreducible, infinite, non-affine Coxeter group is an infinite set.
This implies that an irreducible, infinite Coxeter group is affine if and only if
it contains an abelian subgroup of finite index.

 \vskip 0.3cm

 \vskip
0.3cm

}
\end{minipage}
\end{center}
\vskip 1cm
\baselineskip 16pt
\section{Introduction}
A  \emph{Coxeter system} $(W,S)$ is a
group $W$ and a set $S$ of generators such that
$W$ has a presentation of the form,
\begin{equation}
\label{coxdef}
 W = \langle\ S | (st)^{m_{s t}}=1,\mbox{  }    s, t\in S\rangle,
\end{equation}
where $m_{s t}= m_{t s}$ is a positive integer or $\infty$, and
$m_{s t}=1$ if and only if $s=t$ (A relation $(s t)^{\infty}=1$ is
interpreted as vacuous).  $W$ is called a \emph{Coxeter group}.
The cardinality $|S|$ of $S$ is called
the \emph{rank} of $W$. We are mainly interested in
Coxeter groups of finite rank. So, we assume $|S|$  is finite in this
paper.

For a Coxeter system $(W,S)$, its \emph{Coxeter graph} is a graph
with a vertex set  $S$, and  with two vertices $s\neq t$ joined by
an edge whenever $m_{st}\geq 3$. If $m_{st}\geq 4$, the
corresponding edge is labeled by $m_{st}$. We say that a Coxeter
group $(W,S)$ is \emph{irreducible} if its Coxeter graph is
connected.

Associated to a Coxeter group $(W,S)$, there is a symmetric
bilinear form on a real vector space  $V$, having a basis
$\{\alpha_{s}| s\in S\}$ in one-to-one correspondence with $S$.
The bilinear form ( , ) is defined by setting
\begin{equation}
\label{linearfrom} \mbox{(}\alpha_s, \alpha_t\mbox{)} =
-\cos\frac{\pi}{m_{s t}}.
\end{equation}
The value on the right-hand side is interpreted to be $ -1$ when $m_{s t}=\infty$.

A well-known fact is that a Coxeter group $W$ is finite if and
only if its bilinear form is positive definite. It is stated in
\cite{bourbaki} (page 137) that an irreducible, infinite Coxeter
group has a trivial center, and a proof using the canonical
representations of a Coxeter group, developed by J. Tits (see \cite{bourbaki}
\cite{hump}), is
suggested.

If the bilinear form of an irreducible Coxeter group $(W,S)$ is
positive semi-definite but not positive definite, then $W={\bf \rm
Z}^n\rtimes W_0$, where $W_0$ is a finite Coxeter group and
$n=|S|-1$. We call $W$ an \emph{irreducible, infinite, affine}
Coxeter group in this situation.
  A natural and interesting question, which
was proposed to the author by M. Davis and T. Januszkiewicz, is to
determine if the center of a finite index subgroup of an
\emph{irreducible, infinite, non-affine} Coxeter group is trivial.
By ``\emph{non-affine}''  we mean its bilinear form is neither
positive definite nor positive semi-definite.  The answer is yes.
\begin{theorem}
\label{center} The center of any finite index subgroup of an
irreducible, infinite, non-affine Coxeter group  is trivial.
\end{theorem}

The solution of this question was inspired by a preprint
\cite{paris1} of L. Paris. In \cite{paris1}, by studying the
essential elements (which will be defined in Section 3) of a Coxeter
group, Paris obtained several interesting results on irreducible
Coxeter groups. One of them is  that any irreducible, infinite
Coxeter group cannot be written as a product of two nontrivial
subgroups.

The idea of studying essential elements (Krammer \cite{krammer},
Paris \cite{paris1}) is important in the proof of Theorem 1.1. In
addition, the author makes use of some arguments similar to those in
the proofs of the Flat Torus Theorem and the Solvable Subgroup
Theorem of CAT(0) spaces. For a detailed description of CAT(0)
spaces, the reader is referred to \cite{bridson}. We will explain
briefly in Section 2 a geometric construction   associated to a
Coxeter system $(W,S)$ (see \cite{Davis1}, \cite{Davis2}), which yields a
PE cell complex $\Sigma=\Sigma(W,S)$ (here PE stands for
``piecewise Euclidean''), now commonly called
\emph{Davis complex}. It is proved by G. Moussong \cite{moussong}
that $\Sigma$ is a CAT(0) space and $W$ acts properly and
cocompactly on $\Sigma$ by isometries.

The proof of Theorem 1.1 also relies on the general theory of root
system of a Coxeter group (see Bourbaki  \cite{bourbaki},
 Deodhar \cite{deodhar1} and Krammer
\cite{krammer}). Deodhar \cite{deodhar2} and M. Dyer \cite{dyer1}
independently proved a theorem, which says that any subgroup
generated by a collection of
 reflections in a Coxeter group is
a Coxeter group by itself. This theorem also plays an important role
in the proof.

Using some arguments similar to the proof of Theorem 1.1, we obtain
the following.
\begin{theorem}
\label{product}
Any finite index subgroup of an
irreducible, infinite, non-affine Coxeter group cannot be
expressed as a product of two nontrivial subgroups.
\end{theorem}

After proving these, the author discovered that in a  revised
version of \cite{paris1}, Paris extended his discussions to include
 the
situation of Theorem 1.1 and 1.2 using purely algebraic arguments.
It appears that we both realized the necessity of use of the
reflection subgroup theorem obtained by Deodhar and Dyer to achieve
this aim.  A different proof of Theorem 1.2 can be found in de Cornulier
and de la Harpe \cite{harpec}, where they mention that Theorem 1.1 can be
obtained from a result in Benoist and de la Harpe \cite{harpeb}.

Based on Theorem 1.1 and 1.2,  we obtain the following result.
\begin{theorem}
If  a group $G$ is a direct product of  $n$ irreducible, infinite,
non-affine Coxeter groups, then any finite index subgroup $H$ of $G$
has a trivial center, and $H$ can be expressed uniquely as a direct
product of $m$ nontrivial subgroups of  $H$ (up to the rearrangement
of factors), where each factor cannot be further decomposed and
$1\leq m \leq n$.
\end{theorem}

There are examples, which will be explained in Section 4,
showing that the situation $m<n$  can happen.

Using some arguments in the proof of Theorem 1.1, we obtain another
result, which has an implication on the group ring $R[W]$, where $W$
is an irreducible, infinite, non-affine Coxeter group and $R$ is a
commutative ring with identity. Here $R[W]$ is the free $R$-module
generated by elements of $W$. Any element of $R[W]$ is of the form
$\sum\limits_{w\in W}a(w)w$, where $a(w)\in R$ and $a(w)=0$ for
almost all $w$. The multiplication in $W$ extends uniquely to a
$R$-bilinear product $R[W]\times R[W]\rightarrow R[W]$. This makes
$R[W]$ a ring. First, we have the following.

\begin{theorem}
\label{conjugacyclass} Given an irreducible, infinite, non-affine
Coxeter group $(W,S)$ and any $w\neq 1$ in $W$, the cardinality of
the set $\{gwg^{-1}| g\in W\}$ is infinite.
\end{theorem}

As a comparison, recall that an irreducible, infinite, affine
Coxeter group $W$ has a decomposition $\mathbf{Z}^n\rtimes W_0$,
where $W_0$ is a finite Coxeter group. In this situation,  the
cardinality of the set $\{gwg^{-1}|g\in W\}$ is finite for
$w\in\mathbf{Z}^n$. In summary, Theorem 1.1, 1.2 and 1.4 illustrate the
group theoretic differences between irreducible, infinite,
non-affine Coxeter groups and irreducible, infinite, affine Coxeter
groups, even though the classification between ``affine'' and ``non-affine''
is based on generators and relations and the associated bilinear form. Indeed,
the following fact can be easily proved.

\begin{corollary}
An irreducible, infinite
Coxeter group $W$ is affine if and only if
it contains an abelian subgroup of finite index.
\end{corollary}

\begin{proof}
The ``only if'' part is obvious. For the ``if'' part, suppose that
$A$ is a finite index abelian subgroup of the irreducible,
infinite Coxeter group $W$. Since the number of distinct left
cosets of $A$ in $W$ is the same as the number of distinct right
cosets of $A$ in $W$, we assume that $W=\bigcup\limits_{i=1}^n
w_iA$, where $n$ is a positive integer, $\{w_1, w_2, \cdots ,
w_n\}$ is a designated finite subset of $W$. Given $a\neq 1\in A$,
we have $\{waw^{-1}| w \in W\}=\{w_iaw_i^{-1}| i=1, \cdots, n\}$.
By Theorem 1.4, $W$ cannot be non-affine.  Hence, $W$ is affine.
\end{proof}

The corollary for the group ring $R[W]$ is as follows, the proof of which is left to
the reader.

\begin{corollary}
Let $W$ be an irreducible, infinite, non-affine Coxeter group and
$R$ be a commutative ring with identity. The center of $R[W]$ is
$R$.
\end{corollary}

The paper is organized  as follows. In Section 2,  we state some
basic facts about the combinatorial theory of a Coxeter group. Using
these, the author gives another proof of the statement
 that an irreducible, infinite Coxeter group
has a trivial center. This proof does not use the canonical
representations of Coxeter groups and is of purely combinatorial
nature.
  Then we describe briefly the construction of the Davis complex of a Coxeter
  group.
   The root system of a Coxeter group is introduced in Section 3,
   where some key results from
Krammer's thesis are stated. These are important in our discussions.
The proofs of these theorems are given in Section 4. The last
section is an appendix, where the author supplies some missing
arguments for an important statement (Theorem 3.6) in Section 3.

\section{Basic Combinatorial Theory of Coxeter Groups}
A Coxeter group may be characterized by some combinatorial
conditions, which are stated below. For now, let $W$ be a group
generated by a subset $S$
 of \emph{involutions} (elements of order 2). The \emph{length} $l(w)$ or $l_S(w)$ of an element
$w\in W$,  with respect to $S$,  is the smallest number $d$ such that $w=s_1\cdots s_d$,
with all $s_i\in S$.  This expression for $w$ is called a \emph{reduced decomposition}
 of $w$
 if $d=l(w)$.

 Consider the following conditions.

{\bf (D)  Deletion Condition}. \emph{If $w=s_1\cdots s_d$ with $d>l(w)$, then there are indices
$i<j$ such that $w=s_1\cdots \hat{s}_i\cdots \hat{s}_j\cdots s_d$. where the hats indicate deleted
letters.}

{\bf (E) Exchange Condition}. \emph{Given $w\in W$, $s\in S$, and any reduced decomposition
$w=s_1\cdots s_d$ of $w$, either $l(sw)=d+1$ or else there is an index $i$ such that
$w=ss_1\cdots\hat{s}_i\cdots s_d$.}

{\bf (F) Folding Condition}. \emph{Given $w\in W$ and $s,t\in S$ such that $l(sw)=l(w)+1$ and
$l(wt)=l(w)+1$, either $l(swt)=l(w)+2$ or else $sw=wt$}.

The proof of the following theorem can be found in \cite{bourbaki},
\cite{kbrown} or \cite{davis1}.

\begin{theorem}
\label{coxeter}
A group $W$  generated by a set $S$ of  involutions gives a Coxeter system $(W,S)$ if and only
if $W$ satisfies any one of the  conditions (D), (E) and (F), with the length function
$l(w)=l_S(w)$ defined as above.
\end{theorem}

Given a Coxeter system $(W,S)$, for each subset $T$ of $S$, let
$W_T$ be the subgroup generated by  $T$.  Call it a \emph{special
subgroup} of $W$. Then any element $w\in W$ can be expressed as
$w=w_0 a$ where $a\in W_T$ and $w_0$ is the shortest element in the
left coset $wW_T$.  $w_0$ is characterized by the property
$l(w_0t)=l(w_0)+1$ for any $t\in T$ and is unique in $wW_T$. We say
$w_0$ is \emph{$(\emptyset, T)$-reduced} in this situation. It is clear
this type of  decomposition for $w$ is unique and $w_0$ satisfies
that $l(w_0b)=l(w_0)+l(b)$  for any $b\in W_T$. Similar discussions
for right cosets  give a ``right-hand version'' of the decomposition and
the definition of \emph{$(T, \emptyset)$-reduced} elements.

For $w\in W$, define a subset In$(w)$ of $S$ by
\[ {\rm In}(w)=\{s\in S| l(ws)=l(w)-1\},
\]
and put
\[{\rm Out}(w)=S-{\rm In}(w).\]

We collect some basic facts of finite special subgroups of a
Coxeter group.

\begin{lemma}
\label{long} Suppose $W_T$ is a finite subgroup, where $T\subset
S$. Then there is a unique element $w_T$ in $W_T$ of longest
length. Moreover, the following statements are true.

(1) $w_T$ is an involution.

(2) For any $x\in W_T$, $x=w_T$ if and only if \mbox{ } ${\rm
In}(x)=T$.

(3) For any $x\in W_T$, $l(w_Tx)=l(w_T)-l(x)$.
\end{lemma}

This lemma is taken from exercises in Chapter 4 of \cite{bourbaki}.
The proof of this lemma and the following two lemmas can be found in
Chapter 3 of Davis \cite{davis1}.  For readers' convenience, we
include the proof of Lemma 2.4.

\begin{lemma}
\label{finite} {\rm (\cite{Davis3})} For any $x\in W$, $W_{{\rm
In}(x)}$ is a finite subgroup.
\end{lemma}

\begin{lemma}
\label{davislemma} {\rm (\cite{Davis3})} If $W_T$ is a finite
subgroup of $W$ and $w_T$ is the longest element in $W_T$, then for
$s\in S-T$, $sw_T=w_Ts$ if and only if $m_{st}=2$ for all $t\in T$.
\end{lemma}
\begin{proof}
 If $m_{st}=2$, then $s$ and $t$ commute. Hence, if
$m_{st}=2$ for all $t\in T$, then $s$ and $w_T$ commute.

Conversely, suppose $s$ and $w_T$ commute, where $s\not\in T$. Then
$l(w_Ts)=l(w_T)+1$, so $s\in {\rm In}(w_Ts)$. Since $w_Ts=sw_T$,
$T\subset{\rm In}(w_Ts)$. Therefore, ${\rm In}(w_Ts)=T\cup\{s\}$,
$w_{T\cup\{s\}}=w_Ts$. We want to show that $m_{st}=2$ for all $t\in
T$. Suppose, to the contrary, that $m_{st}>2$, for some $t\in T$.
Then $l(sts)=3$, $l((w_Ts)(sts))=l(w_T)+1-3=l(w_T)-2$ by Lemma
~\ref{long}. On the other hand,
$l((w_Ts)(sts))=l(w_Tts)=l(w_Tt)+1=l(w_T)$, a contradiction. Hence,
the conclusion of Lemma 2.4 holds.

%Since $\{s,t\}\subset {\rm In}(w_Ts)$, which generates a finite
%subgroup (Lemma ~\ref{finite}), $m_{st}\neq\infty$. Consider the
%element $sts$. Since $m_{st}>2$,
\end{proof}

Now we prove the claim.
\begin{proposition}
\label{trivialcen} The center of an irreducible, infinite Coxeter
group $(W,S)$ is trivial.
\end{proposition}

\begin{proof}
If $w\neq 1$ is in the center of $W$, then $ws=sw$ for any $s\in S$.
Put $S_1={\rm In}(w)$ and $S_2={\rm Out}(w)$.  Then
$S_1\neq\emptyset$. Write $w=w_0a$ with $a\in W_{S_1}$ and $w_0$
being $(\emptyset, S_1)$ reduced. Notice that
$l(w)-1=l(ws)=l(w_0as)$ for any $s\in S_1$, it follows that
$l(as)=l(a)-1$ for all $s\in S_1$.  By Lemma 2.2 and Lemma 2.3, $a$
is the (unique) longest element in the finite subgroup $W_{S_1}$,
and $a^2=1$.

Now, continue our discussions and consider the ``right-hand version''
of the above-mentioned decomposition of $w$. Since $w$ is in the
center of $W$, we have $w=aw_1$, where $w_1$ is $(S_1,\emptyset)$
reduced, $a$ is the longest element in $W_{S_1}$.  Hence,
$w_0=wa=aw=w_1$, $w=aw_0=w_0a$.

Notice that for any $t\in S_2$, $l(wt)=l(w)+1$, it follows that
$l(w_0t)=l(w_0)+1$, since otherwise we would have $l(wt)\leq
l(a)+l(w_0t)\leq l(w)-1$, contradicting the definition of $S_2$.
Therefore,  $w_0$ is $(\emptyset, S_2)$-reduced, and hence is
$(\emptyset, S)$-reduced.  This implies that $w_0=1$ and $w=a$,
i.e.,
 $w=w_{S_1}$. So, $w_{S_1}$ commutes with every element in
$S_2=S-S_1$.  Now, Lemma \ref{davislemma} implies that  $m_{st}=2$
for any $s\in S_1$, $t\in S_2$. The irreducibility of $W$ implies
$S_2=\emptyset$ and hence $W$ is a finite Coxeter group, a
contradiction. This finishes the proof of Proposition
~\ref{trivialcen}.
\end{proof}

To prove the theorems stated in the introduction, we need the fact that a
Coxeter group acts properly and cocompactly on a CAT(0) space. Here
we give a brief description of the Davis complex. For  a more
complete account of it, the reader is referred to \cite{Davis1}
\cite{Davis2} \cite{Davis3}.

Let $(W,S)$ be a Coxeter system. We define a poset, denoted  ${\rm
S}^f(W,S)$ (or simply ${\rm S}^f$ ), by putting
\[{\rm S}^f=\{T | T\subset S \mbox{ and } W_T \mbox{ is finite}\}.\]
This poset is partially ordered by inclusion. It is clear that ${\rm S}^f-\{\emptyset\}$
is isomorphic to the poset of simplices of an abstract simplicial
complex, which is denoted by $N(W,S)$ (or simply $N$). $N$ is called
the \emph{nerve} of $(W, S)$.
\begin{theorem}
\label{moussong} $(${\rm Gromov, Moussong \cite{Davis2}
\cite{moussong}}$)$. Associated to a Coxeter system $(W,S)$, there
is a PE cell complex $\Sigma(W,S)\mbox{ } (=\Sigma)$ with the
following properties.

(1) The poset of cells in $\Sigma$ is the poset of cosets
\[W{\rm S}^f=\coprod\limits_{T\in{\rm S}^f}W/W_T.\]

(2) $W$ acts by isometries on $\Sigma$ with finite stabilizers and with
compact quotient.

(3) Each cell in $\Sigma$ is simple (so that the link of each vertex is a simplicial
cell complex).  In fact, this complex is just $N(W,S)$.

(4) $\Sigma$ is CAT(0).
\end{theorem}

\section{Root System and Essential Elements of a Coxeter Group}
Recall from the introduction that  for a Coxeter system $(W,S)$,
there is a symmetric bilinear form ( , ) on a real vector space $V$,
having a basis $\{\alpha_{s}| s\in S\}$ in one-to-one correspondence
with $S$.

Now, for each $s\in S$,  define a linear transformation $\sigma_s:
V\rightarrow V$ by $\sigma_s
\lambda=\lambda-2\mbox{(}\alpha_s,\lambda\mbox{)}\alpha_s$. Then
$\sigma_s$ is a linear reflection. It has order 2 and fixes the
hyperplane $H_s=\{\delta\in V|\mbox{(}\delta,\alpha_s\mbox{)}=0\}$
pointwise,  and $\sigma_s\alpha_s=-\alpha_s$. We have the following
theorem (see \cite{deodhar1}, \cite{hump}).

\begin{theorem}
\label{repre1}
There is a unique homomorphism $\sigma : W \rightarrow \mbox{GL(}V\mbox{)}$
sending $s$ to $\sigma_s$. This homomorphism is a faithful representation of
$W$ and the group $\sigma\mbox{(}W\mbox{)}$ preserves the bilinear form.
 Moreover, for each pair $s, t\in S$, the order of $st$ in $W$ is
precisely $m_{s t}$.
\end{theorem}

From now on, we write $w(\alpha)$ for $\sigma(w)(\alpha)$, when $\alpha\in V$
and $w\in W$.

 The \emph{root system}  $\Phi$ of $W$, is defined to be the
collection of all vectors $w\mbox{(}\alpha_s\mbox{)}$, where $w\in W$ and $s\in S$.
An important fact about the root system is that any root $\alpha\in \Phi$ can be
expressed as
\[\alpha=\sum\limits_{s\in S} c_s \alpha_s, \]
where all the coefficients  $c_{s}\geq 0$ (we call $\alpha$ a
positive root and write $\alpha>0$), or all the
 coefficients $c_{s}\leq 0$ (call $\alpha$ a negative root and write $\alpha <0$).
 Write $\Phi^{+}$ and $\Phi^{-}$ for the respective sets of positive and negative roots.
 Then $\Phi^{+}\bigcap\Phi^{-}=\emptyset$ , $\Phi^{+}\bigcup\Phi^{-}=\Phi$
 and $\Phi^{-}=-\Phi^{+}$. The map from $\Phi$ to $R=\{wtw^{-1}| w\in W, t\in S\}$
 (the set of reflections in $W$) given
 by $\alpha=w(\alpha_s)\mapsto  wsw^{-1}$ is well-defined
and restricts to a bijection from $\Phi^{+}$ ($\Phi^{-}$) to $R$, and $\sigma(wsw^{-1})=t_\alpha$,
where $t_\alpha$ is the linear reflection given by
 $t_{\alpha}\lambda=\lambda-2(\alpha,\lambda)\alpha$. The following
fact is important when discussing root systems.

\begin{proposition}
\label{length} {\rm (\cite{deodhar1}, \cite{hump})}   Let $w\in W$,
$\alpha\in\Phi^{+}$. Then $l(wt_{\alpha})>l(w)$ if and only if
$w(\alpha)>0$.
\end{proposition}

 The statements in the remaining part of this section are mostly due to Krammer
 \cite{krammer}, as revised by Paris \cite{paris1}. Let $u$, $v\in W$
 and $\alpha\in \Phi$. We say that $\alpha$ \emph{separates} $u$ and
 $v$ if $u\alpha\in\Phi^{\epsilon}$ and
 $v\alpha\in\Phi^{-\epsilon}$, where $\epsilon\in\{+,-\}$. Let $w\in
 W$ and $\alpha\in\Phi$. We say that $\alpha$ is $w$-\emph{periodic} if
 there is some positive integer $m$ such that $w^m\alpha=\alpha$.

 \begin{lemma}
 {\rm (\cite{paris1})}   Let $w\in W$ and $\alpha\in \Phi$. Then exactly
 one of the following holds.

 (1) $\alpha$ is $w$-periodic.

 (2) $\alpha$ is not $w$-periodic, and the set $\{m\in\mathbf{Z}|
 \alpha\mbox{ separates } w^m\mbox{ and } w^{m+1}\}$ is finite and
 has an even cardinality.

 (3) $\alpha$ is not $w$-periodic, and the set $\{m\in\mathbf{Z}|
 \alpha\mbox{ separates } w^m\mbox{ and } w^{m+1}\}$ is finite and
 has an odd cardinality.

 \end{lemma}

 We say that $\alpha$ is $w$-\emph{even} in Case 2, and
 $w$-\emph{odd} in Case 3.

 \begin{lemma}
 {\rm (\cite{paris1})}   Let $\alpha\in\Phi$, $w\in W$, and
 $p\in\mathbf{N}$, $p\geq 1$. Then

 (1) $\alpha$ is $w$-periodic if and only if $\alpha$ is
 $w^p$-periodic.

 (2) $\alpha$ is $w$-even (resp., $w$-odd) if and only if $\alpha$ is
 $w^p$-even (resp., $w$-odd).
 \end{lemma}

A subgroup $G$ is called a \emph{parabolic subgroup} of a Coxeter
group $(W,S)$ if $G=xW_Ix^{-1}$ for some $x\in W$ and $I\subset S$.
An element $w\in W$ is called an \emph{essential element} if it does
not lie in any proper parabolic subgroup of $W$.

Following Krammer \cite{krammer}, we define the \emph{parabolic
closure} ${\rm Pc}(A)$ of a subset $A$ of a Coxeter group $W$ to be
the intersection of all parabolic subgroups containing $A$. In
\cite{qi1}, the author gives a proof of the conclusion that ${\rm
Pc}(A)$ is a parabolic subgroup of $W$. Using this terminology,  an
element $w\in W$ is essential if and only if ${\rm Pc}(w)=W$.

 Paris shows the existence of essential elements in \cite{paris1}.
\begin{proposition}
\label{essentialele} {\rm (Paris \cite{paris1})}   Given a Coxeter
group $(W,S)$, where $S=\{s_1, s_2, \cdots, s_n\}$, then
$c=s_n\cdots s_2s_1$ (which is called a Coxeter element) is an
essential element of $W$.
\end{proposition}

In \cite{paris1} Paris attributes the following  result
to Krammer  \cite{krammer}.
\begin{theorem}
\label{krammeressl} For an irreducible, infinite Coxeter group
$(W,S)$, an element $w\in W$ is essential if and only if $W$ is
generated by the set $\{t_{\alpha}| \alpha\in\Phi^{+} \mbox{ and }
\alpha\mbox{ is w-odd}\}$.
\end{theorem}

To the author's understanding, this statement does not appear in
Krammer's thesis \cite{krammer}, but it can be proved using some
results that Krammer has established. The author will give the proof
of this theorem in the Appendix.

The next result is obvious from the discussions in the preceding
paragraphs.
\begin{corollary}
\label{essentialelem}{\rm (\cite{paris1})}  Assume that $W$ is an
irreducible, infinite Coxeter group. Let $w\in W$ and $p$ be a
positive integer. Then $w$ is essential if and only if $w^p$ is
essential.
\end{corollary}

The following theorem, which appears in Krammer \cite{krammer}
(page 69), is very important to us.
\begin{theorem}
Assume that  $W$ is an irreducible, infinite, non-affine Coxeter
group. Let $w\in W$ be an essential element. Then $\langle
w\rangle=\{w^m| m\in\mathbf{Z}\}$ is a finite index subgroup of the
centralizer $C(w)$ of $w$ in $W$.
\end{theorem}

Now we come to the proofs of the theorems stated in the
Introduction.
\section{Proofs of the Theorems}
We begin with a lemma.
\begin{lemma}
\label{dinfty} If an irreducible Coxeter group $(W, S)$ contains
an infinite cyclic subgroup $\langle x\rangle$ as a finite index
subgroup, then $W\cong D_\infty$, the infinite dihedral group.
\end{lemma}
\begin{proof}
The key point is that $W$ acts on the CAT(0) space
$\Sigma=\Sigma(W,S)$ properly and cocompactly. What follows is
similar to the proof of the Flat Torus Theorem for CAT(0) spaces
(see page 246 in \cite{bridson}). Since $[W:\langle x\rangle]$ is
finite, there is a positive integer $k$ such that $\langle
x^k\rangle$ is a normal subgroup.  It is a finite
index subgroup of $W$. ${\rm \bf Min}(\langle x^k\rangle)$ is a
non-empty
 closed  subspace of $\Sigma$ isometric to one of the form
 $Y\times\mathbf{R}$, where $Y$ is a closed convex subspace of $\Sigma$,
 $\mathbf{R}$ is
 the set of real numbers
  (For any group
$\Gamma$ acting on a CAT(0) space by isometries, ${\rm \bf
Min}(\Gamma)$ is a certain subspace as defined on page 229 in
\cite{bridson}).  $x^k$ acts as a nontrivial translation on
 the factor $\mathbf{R}$ and acts
as an identity map on $Y$. Since $\langle x^k\rangle\subset W$ is
normal, $W$ acts by isometries of ${\rm \bf Min}(\langle
x^k\rangle)$, preserving the splitting. By properties of CAT(0)
spaces, the fixed point set $Y_1$ of the induced action of the
finite group $W/\langle x^k\rangle$ on the complete CAT(0) space
$Y$ is a non-empty, closed, convex subset of $Y$. By construction
$Y_1\times\mathbf{R}$ is $W$-invariant and the action of $W$ on
the factor
$Y_1$ is the identity. Pick $y\in Y_1$ and consider the action of
$W$ on $L=\{y\}\times\mathbf{R}$. The restriction of each $s$ to
$L$ is either a reflection
 or the identity. Since $x^k$ is a translation on $L$, there are elements
  $s_1$ and $s_2$
 in $S$, which act
 as different reflections on $L$. So, the order of $s_1s_2$ is infinite. By
  considering the collection of cosets
  $\{\langle x^k\rangle (s_1s_2)^l|l\in\mathbf{Z}\}\subset W/{\langle x^k\rangle}$,
   we know there is a positive integer $d$ such that $(s_1s_2)^d\in \langle x^k\rangle$.
   Hence,  $[W:\langle s_1s_2\rangle]$  is finite,  and therefore, so is $[W:W_{\{s_1, s_2\}}]$.
   This is
   impossible unless $W=W_{\{s_1, s_2\}}=D_{\infty}$, because of the following
   result
   proved by Deodhar \cite{deodhar1}: \emph{For an
   irreducible, infinite Coxeter group $(W,S)$ and any proper subset $J$ of $S$,
   $[W:W_J]$ is infinite} (This Deodhar's Theorem was discovered again in T. Hosaka   \cite{hosaka1},
   using different methods).
\end{proof}

\noindent \emph{Proof of Theorem 1.1}. Assume that $(W,S)$ is an
irreducible, infinite, non-affine Coxeter group and $G$ is a finite
index subgroup of $W$. Let $Z(G)$ be the center of $G$. The proof is
divided into two cases.

Case 1. $Z(G)$ contains an element $x$ of infinite order. By
Proposition ~\ref{essentialele}, we can pick an essential element
$w$ of $W$. Since the number of cosets $\{Gw^m|m\in\mathbf{Z}\}$
is finite, there is a positive integer $p$ such that $w^p\in
G$. By Corollary ~\ref{essentialelem}, $w^p$ is essential in $W$.
Since $x$ is in the centralizer $C(w^p)$ of $w^p$, the number of
cosets $\{\langle w^p\rangle x^m|m\in\mathbf{Z}\}\subset
C(w^p)/{\langle w^p\rangle}$ is finite, so, there is a positive
integer $q$ such that $x^q\in\langle w^p\rangle$. Now, $x^q$ and
hence,  $x$ is essential in $W$. Since $G$ is contained in the
centralizer $C(x)$ of $x$, it follows that the index $[G:\langle
x\rangle]\leq[C(x):\langle x\rangle]$ is finite  by Theorem 3.8.
Therefore, $[W:\langle x\rangle]$ is finite. Now, Lemma
~\ref{dinfty} implies that $W$ is the infinite dihedral group
$D_{\infty}$, which is impossible since $D_\infty$ is affine.

Case 2. $Z(G)$ is a torsion subgroup, i.e., every element in
$Z(G)$ has a finite order. A preliminary result used in the proof
of the Solvable Subgroup Theorem for CAT(0) spaces (see page 247
in \cite{bridson}) states that \emph{if a group $\Gamma$ acts
properly and cocompactly by isometries on a CAT(0) space, then
every abelian subgroup of $\Gamma$ is finitely generated}. It
follows  that $Z(G)$ is finitely generated and hence, is finite. A
result of Tits states that \emph{any finite subgroup of a Coxeter
group is contained in a finite parabolic subgroup}. Therefore, in
this case, the parabolic closure ${\rm Pc}(Z(G))$ of $Z(G)$ is a
finite parabolic subgroup of $W$. Without loss of generality, we
may assume that ${\rm Pc}(Z(G))$ is a finite special parabolic
subgroup $W_K$, where $K\subset S$.

Since $Z(G)$ is normal in $G$, $gW_K g^{-1}$ is a parabolic
subgroup containing $Z(G)$ for any $g\in G$. By the uniqueness of
the parabolic closure (or by the discussion of the rank of the
intersection of two parabolic subgroups in \cite{qi1}), we have
$gW_K g^{-1}=W_K$ and hence, $G\subset N(W_K)$ (the normalizer of
$W_K$ in $W$). Therefore, $[W:N(W_K)]$ is finite. This implies
that the set $R_1=\{wtw^{-1}|t\in K, w\in W\}$ is finite. Now,
consider the reflection subgroup $W_1$ of $W$ generated by $R_1$.
$W_1$ is a Coxeter group by \cite{deodhar2} or \cite{dyer1}, with
a set $S_1$ of distinguished generators, where
$S_1\subset\bigcup\limits_{w_1\in W_1}w_1R_1{w_1}^{-1}$ (It is
clear that $\bigcup\limits_{w_1\in W_1}w_1R_1{w_1}^{-1}= R_1$ in
the present situation). Hence, the set of reflections in $W_1$,
which by definition is $\{w_1t_1w_1^{-1}|t_1\in S_1, w_1\in
W_1\}(\subset R_1)$, is finite. Therefore, $W_1$ is a finite
Coxeter group. Suppose that ${\rm Pc}(W_1)=yW_Ly^{-1}$, where
$y\in W$ and $L$ is a proper subset of $S$. Since $W_1$ is normal
in $W$, $yW_Ly^{-1}$ is a (proper) normal subgroup of $W$. This is
impossible, due to the  result proved by Paris \cite{paris1}:
\emph{any proper nontrivial special subgroup of an irreducible
Coxeter group is not normal}. In the current situation, we need to
replace the distinguished set $S$ of generators by $ySy^{-1}$.

In conclusion, $Z(G)=\{1\}$. This finishes the proof of Theorem
1.1.

\begin{Remark} Paris obtains the above mentioned result based on R.
Howlett's description of the normalizer of a special subgroup of a
Coxeter group in \cite{howlett1}, where Howlett discusses the
situation of finite Coxeter groups.  In fact, Howlett's description is
 valid in general. The proof is just a slight modification of
that for finite Coxeter groups.
\end{Remark}

The proof of Theorem 1.2 is similar.

\noindent \emph{Proof of Theorem 1.2}. Let $G$ be a finite index
subgroup of an irreducible, infinite, non-affine Coxeter group.
Suppose that $G=A\times B$, where $A$
and $B$ are nontrivial subgroups. Pick an essential element $w$ in
$W$. Following the idea in the proof of Theorem 1.1, we know that
there is a positive integer $p$ such that $w^p\in G$. Now $w^p$ is
essential and $w^p=ab$ for some $a\in A$ and $b\in B$. At least
one of them, say, $a$ has an infinite order, because $w^p$ has an
infinite order, and $a$ and $b$ commute. Notice that since
$aw^p=w^pa$, $a$ is in the centralizer $C(w^p)$ of $w^p$. By
considering the collection of cosets $\{\langle w^p\rangle
a^m|m\in\mathbf{Z}\}$, we conclude that there is a positive
integer $q$ such that $a^q\in \langle w^p\rangle$. So, $a$ is an
essential element. Since each element of $B$ commutes with $a$,
$B\subset C(a)$, the centralizer of $a$. Since the collection of
cosets $\{\langle a\rangle h|h\in B\}$ is finite and $A\cap
B=\{1\}$, we know that $B$ is finite. Without loss of generality,
we can assume that ${\rm Pc}(B)=W_I$,  a finite parabolic
subgroup. Since $B$ is normal in $G$, the uniqueness of parabolic
closure implies that $G\subset N(W_I)$. So, $[W:N(W_I)]$ is
finite. The arguments in Case 2 of the proof of Theorem 1.1 now
apply to derive a contradiction. Hence, the conclusion of Theorem
1.2 holds. \vskip 0.4cm

\noindent \emph{Proof of Theorem 1.3}. Assume that
\begin{equation}
\label{decomp1} G=W_1\times W_2\times \cdots\times W_n\mbox{
(internal direct product)},
\end{equation}
 where each $W_i$ is an irreducible, infinite,
non-affine Coxeter group. Let $H$ be a finite index subgroup of
$G$. Denote by $p_i$ the projection $p_i: G \rightarrow
W_i$ of $G$ to the $i$th factor $W_i$. Let
$c=\prod\limits_{i=1}^n c_i\in Z(H)$, the center of $H$, where
$c_i\in p_i(H)\subset W_i$. Pick an arbitrary
$h=\prod\limits_{i=1}^n h_i\in H$ with $h_i\in p_i(H)$. Then
$ch=hc$. This implies that $\prod\limits_{i=1}^n(c_i
h_i)=\prod\limits_{i=1}^n(h_i c_i)$, so, $c_i h_i= h_i c_i$ for
any $i$. Notice that since $h_i\in p_i(H)$ is arbitrary,  $c_i\in
Z(p_i(H))$. Since $W_i=p_i(G)$ and $[p_i(G): p_i(H)] =[G: p_i^{-1}(p_i(H))]\leq [G:
H]<\infty$, Theorem 1.1 implies that $c_i=1$. Hence, $Z(H)=\{1\}$.

Now we use induction on $n$, the number of factors in expression
(\ref{decomp1}), to prove that if a finite index subgroup $H$ of
$G$ can be expressed as
\begin{equation}
\label{decomp2} H=H_1\times H_2\times\cdots\times H_m\mbox{
(internal direct product)},
\end{equation}
where each $H_i$ is a nontrivial subgroup, then $m\leq n$.

The case $n=1$ is just Theorem 1.2. Now, assume $n\geq 2$ and
$m\geq 2$. Notice that $p_1(H)=p_1(H_1)p_1(H_2)\cdots p_1(H_m)$, and
\[p_1(H_i)\bigcap \{p_1(H_1)\cdots p_1(H_{i-1})p_1(H_{i+1})\cdots
p_1(H_m)\}=\{1\},\]
 because this intersection is contained in
$Z(p_1(H))$, which is trivial by Theorem 1.1 (knowing that
$p_1(H)$ is a finite index subgroup of $W_1$). Hence
\[p_1(H)=p_1(H_1)\times\cdots\times p_1(H_m).\]
By Theorem 1.2 only one of the factors on the right-hand side,
say, $p_1(H_1)$ can be nontrivial, and all other $p_1(H_j)$
($j\neq 1$) are trivial. So,  $p_1(H_1)=p_1(H)$.  Without loss of
generality, we can assume that $p_i(H_1)$ is nontrivial for $i=1,
\ldots, l$, and is trivial for $i=l+1, \ldots, m$. This implies
that $p_i(H_j)=\{1\}$ for $i=1, \ldots, l$, $j\not=1$. Hence,
\[H_1\subset W_1\times\cdots\times W_l,\mbox{ }
H_2\times\cdots\times H_m\subset W_{l+1}\times\cdots\times W_n.\]

Now use the induction hypothesis and the following simple fact
(the proof of which is left to the reader),
\begin{lemma}
Let $G_1$ and $G_2$ be two groups. If $N_i$ is a subgroup of
$G_i$, $i=1, 2$ and $[G_1\times G_2: N_1\times N_2]<\infty$, then
$[G_i: N_i]<\infty$, $i=1, 2$, and $[G_1\times G_2: N_1\times
N_2]=[G_1: N_1] [G_2: N_2]$.
\end{lemma}
 We conclude that $m\leq n$.

Having this inequality in mind, we may  continue to decompose some
factors in expression (\ref{decomp2}), until each factor cannot be
further decomposed (It is a finite step procedure due to the above
inequality). From now on, when we talk about  a decomposition of
form (\ref{decomp2}), we assume that each factor cannot be further
decomposed.

Suppose that there is another decomposition
\begin{equation}
\label{decomp3} H=K_1\times\cdots\times K_r,
\end{equation}
where each factor $K_j$ cannot be further decomposed. Let $q_j:
H\rightarrow K_j$ be the projection of $H$ onto its
$j$th factor in the decomposition (\ref{decomp3}). We know that
$K_i=q_i(H)=q_i(H_1)q_i(H_2)\cdots q_i(H_m)$. And, notice that
$q_i(H_k)\bigcap\{\prod\limits_{j\neq k} q_i(H_j)\}=\{1\}$,
because this intersection is contained in the center  $Z(K_i)$ of
$K_i$, and $Z(K_i)\subset Z(H)$, while the latter is trivial by
the first part of this theorem, we have the following,
\[K_i=q_i(H_1)\times\cdots\times q_i(H_m).\]
By the assumption that $K_i$ cannot be further decomposed,
$K_i=q_i(H_{\phi(i)})$ for some $\phi(i)$, and $q_i(H_j)$ is trivial
for $j\neq\phi(i)$.  So, $\phi$ defines a map from $\{1, 2, \ldots,
r\}$ to $\{1, 2, \ldots, m\}$, and it is surjective because, for any
$j\in\{1, 2, \ldots, m\}$, there is an $i\in\{1, 2, \ldots, r\}$
such that the restriction $q_i|_{H_j}$ is nontrivial. Hence $m\leq
r$. Similar discussions also yield $r\leq m$. Therefore, $m=r$ and
$\phi$ is a bijection. After re-indexing, we may assume $\phi={\rm
id}$ (the identity map). This means that $H_i\subset K_i$ and the
restriction $q_i|_{H_i}$ is indeed the inclusion $H_i\hookrightarrow
K_i$. Since
\[H=H_1\times\cdots\times H_m=K_1\times\cdots\times K_m,\]
we know that $H_i=K_i$, for $i=1, \ldots, m$. This is the claimed
unique decomposition of $H$ and $m$ is determined uniquely by $H$.
The proof of Theorem 1.3 is completed. \vskip 0.4cm

The situation that $m<n$ may happen. To illustrate this,
we need the following lemma.
\begin{lemma}
Let $W=W_1\times W_2$, where $W_1$ and $W_2$ are irreducible,
infinite, non-affine Coxeter groups. If $H$ is a finite index
subgroup of $W$ and $H=H_1\times H_2$, where $H_i\neq\{1\}$,
$i=1, 2$, then after re-indexing, $H_i\subset W_i$.
\end{lemma}
\begin{proof}
As in the proof of Theorem 1.3, let $p_i$ be the
projection $p_i: W\rightarrow W_i$. It follows that $p_1(H)$ is a
finite index subgroup of $W_1$ and $p_1(H)=p_1(H_1)\times
p_1(H_2)$ by repeating the arguments in the proof of Theorem 1.3.
In this product decomposition, only one factor, say, $p_1(H_1)$ is
nontrivial by Theorem 1.2. So, $p_1(H)=p_1(H_1)$ and
$p_1(H_2)=\{1\}$. This implies that $H_2\subset W_2$. So,
$p_2(H_2)\neq\{1\}$. By Theorem 1.2, $p_2(H_1)=\{1\}$,
$H_1\subset W_1$.
\end{proof}

In the following two examples,
let $(W_i, S_i)$ be an
irreducible, infinite, non-affine Coxeter group, $i=1, 2$. Suppose
that $W=W_1\times W_2$ and $S=S_1\cup S_2$.
\begin{example} (M. Davis) Suppose that $f_i: W_i\rightarrow G$ ($i=1,2$)
is a surjective homomorphism from $W_i$ to a finite group $G$.  Let
$H=\{(w_1, w_2)\in W|f_1(w_1)=f_2(w_2)\}$. Then $[W: H]=|G|$,
$p_1(H)=W_1$ and  $p_2(H)=W_2$, where $p_i$ is the projection defined in
Lemma 4.3. It follows from Lemma 4.3 that $H$  cannot be  further decomposed.
\end{example}
\begin{example}
Let $\phi: W\rightarrow D_1$ be the surjective homomorphism such
that $\phi(s)=-1$ for any $s\in S$, where $D_1=\{-1, 1\}$. Let
$H=\ker(\phi)$. Then $[W:H]=2$. If $H$ had a  decomposition
$H=H_1\times H_2$,  $H_i\neq\{1\}$, then one of the factors, say,
$H_1$ would be $W_1$ because of Lemma 4.3 and Lemma 4.2 and  $[W:
H]=2$. This would imply that $\phi(s_1)=1$ for any $s_1\in S_1$, a
contradiction. Hence, $H=\ker(\phi)$ cannot be further decomposed.
\end{example}

\noindent \emph{Proof of Theorem 1.4}. Let $w\neq 1$ be an element
of $W$. We want to prove that $[W: C(w)]=\infty$. The proof is
divided into two cases.

Case 1. The order of $w$ is finite. In this case, the parabolic
closure ${\rm Pc}(\langle w\rangle)$  of $ \langle w\rangle$ is a
finite parabolic subgroup. Without loss of generality, assume that
${\rm Pc}(\langle w\rangle)=W_K$, where $K\subset S$. The
uniqueness of  the parabolic closure and the fact that
$gwg^{-1}=w$ for any $g\in C(w)$ imply that $gW_K g^{-1}=W_K$ for
$g\in C(w)$. Hence, $C(w)\subset N(W_K)$. The discussion of Case 2
in the proof of Theorem 1.1 yields that $[W: N(W_K)]=\infty$ when
$W$ is an irreducible, infinite Coxeter group and $W_K$ is a
nontrivial finite special subgroup. Therefore, in this situation,
$[W: C(w)]=\infty$.

Case 2. The order of $w$ is infinite. Suppose that $[W:
C(w)]<\infty$. Pick an essential element $x\in W$. Since the number
of cosets $\{C(w)x^k|k\in \mathbf{Z}\}$ is finite, there is a
positive integer $m$ such that $x^m\in C(w)$. By Corollary
\ref{essentialelem}, $x^m$ is essential. Notice that $w\in C(x^m)$
and the number of cosets $\{\langle x^m\rangle w^l|l\in
\mathbf{Z}\}$ in $C(x^m)$ is finite because of Theorem 3.8, we
conclude that there is a positive integer $n$ such that
$w^n\in\langle x^m\rangle$. Now, Corollary \ref{essentialelem}
implies that $w^n$ is essential, and hence, so is $w$. Then, by
Theorem 3.8, we have $[W: \langle w\rangle]<\infty$. It follows from
Lemma \ref{dinfty} that $W=D_\infty$, contradicting the assumption
that $W$ is non-affine.

In conclusion, $[W:C(w)]=\infty$ for $w\neq 1$. The conclusion of
the theorem follows immediately.

\section{Appendix}
In this section we summarize some results from Krammer
\cite{krammer} and supply the missing arguments for the proof of
Theorem 3.6. Recall that we assume that the Coxeter system $(W,S)$
has a finite rank. It has a faithful representation as described in
Section 3.

A $w$-periodic root $\alpha$ is called $w$-\emph{critical}
(\cite{krammer} page 53) if it satisfies the condition: The bilinear
form (defined in the Introduction) is positive definite on
Span$\{\langle w\rangle\alpha\}$ and $\Sigma_{n=1}^kw^n\alpha=0$,
where $k$ is the smallest positive integer with $w^k\alpha=\alpha$.

We need the following results.

\begin{lemma}
\label{criticalr} {\rm (\cite{krammer} page 54)} The subgroup
generated by $\{t_\alpha| \alpha \mbox{ is w-critical}\}$ is a
finite subgroup of $W$.
\end{lemma}

\begin{theorem}
\label{krammerpara} {\rm (\cite{krammer} page 56)} The parabolic
closure {\rm Pc}$(w)$ of $w$ equals the subgroup  of $W$ generated
by the reflections $\{t_\alpha|\alpha \mbox{ is  w-critical or
w-odd}\}$.
\end{theorem}

Let $W_1$ be the subgroup generated by $\{t_\alpha|\alpha \mbox{ is
\emph{w}-odd}\}$ and $\Phi_1=\{w_1(\alpha)|w_1\in W_1, \mbox{
}\alpha \mbox{ is  \emph{w}-}\\
\emph{odd}\}$.

\begin{lemma}
\label{criticalroot} {\rm (\cite{krammer} page 56)} Let $\alpha$
be $w$-periodic, then either $(\alpha,\beta)=0$ for any
$\beta\in\Phi_1$, or $\alpha\in\Phi_1$.
\end{lemma}

\begin{lemma}
\label{radical} {\rm (\cite{hump} page 131)} Assume that $(W,S)$
is irreducible, then any proper $W$-invariant subspace is
contained in the radical $V^{\perp}$ (of the bilinear form),
i.e., $V^{\perp}=\{v\in V|(v,\alpha_s)=0 \mbox{ for any } s\in
S\}$.
\end{lemma}

\noindent \emph{Proof of Theorem ~\ref{krammeressl}}. Proposition
~\ref{essentialele} states that essential elements exist in a
Coxeter group. Now, assume that $(W, S)$ is an irreducible,
infinite Coxeter group. For any essential element $w$ of $W$,
Pc$(w)=W$. It follows from Lemma ~\ref{criticalr} and Theorem
~\ref{krammerpara} that $\Phi_1$ is nonempty. Let $V_0$ be the
subspace spanned by those $w$-critical roots $\alpha$ satisfying
the condition that $(\alpha, \beta)=0$ for any $\beta\in \Phi_1$.
Then Theorem ~\ref{krammerpara} and Lemma ~\ref{criticalroot}
 imply that $V_0$ is a
$W$-invariant subspace. Obviously, $V_0\neq V$, since otherwise,
$(\beta, \beta)=0$ for $\beta\in\Phi_1$,  which is absurd.  It follows from
Lemma ~\ref{radical}  that $V_0\subset V^{\perp}$. This is
impossible unless $V_0=0$, since any critical root $\alpha\in V_0$
satisfies $(\alpha, \alpha)=1$. Therefore,
 if $w\in W$ is an essential element,  $W={\rm Pc}(w)$ is generated by $w$-odd
reflections.  The converse is clear from Theorem
~\ref{krammerpara}.

\vskip 0.6cm

\noindent {\large {\bf Acknowledgements}}

\noindent The author would like to thank Michael Davis and Tadeusz
Januszkiewicz for the inspiration and numerous discussions.

\end{document}